\documentclass[11pt,a4paper]{amsart}
\usepackage[a4paper,margin=1in]{geometry}
\usepackage{amsmath,amssymb,mathtools}
\usepackage{stix2}
\usepackage{graphicx,booktabs,array,longtable,tikz}
\usepackage{microtype}
\usepackage[numbers,sort&compress]{natbib}
\usepackage[colorlinks=true,linkcolor=black,citecolor=black,urlcolor=black,
 pdfauthor={Aryaman Chandra},pdftitle={Near-Return Arithmetic Landscapes of the Discrete Cat Map},
 pdfsubject={Mathematics}]{hyperref}
\usetikzlibrary{arrows.meta,positioning}
\theoremstyle{plain}
\newtheorem{theorem}{Theorem}[section]
\newtheorem{lemma}[theorem]{Lemma}
\newtheorem{proposition}[theorem]{Proposition}
\newtheorem{corollary}[theorem]{Corollary}
\theoremstyle{definition}
\newtheorem{definition}[theorem]{Definition}
\newcommand{\Z}{\mathbb Z}
\newcommand{\R}{\mathcal R}
\newcommand{\Lscape}{\mathcal L}
\newcommand{\Vis}{\operatorname{Vis}}
\newcommand{\ord}{\operatorname{ord}}
\newcommand{\Prob}{\mathbb P}
\newcommand{\E}{\mathbb E}
\numberwithin{equation}{section}
\title[Near returns of the discrete cat map]{NEAR-RETURN ARITHMETIC LANDSCAPES OF THE DISCRETE CAT MAP}
\author{Aryaman Chandra}
\address{Heritage Xperiential Learning School, Gurugram, Haryana, India}
\date{}
\begin{document}
\begin{abstract}
We extend the arithmetic landscape function of a finite dynamical system by incorporating visits to neighbouring points along an orbit. For the discrete cat map, these visits define a near-return graph whose directed lags determine the associated landscape. Under king adjacency on a torus of side $m\geq7$, we prove that every connected component is an isolated vertex, an edge, a two-edge path or a triangle with a leaf. The last occurs precisely when $5\nmid m$ and some power of the cat map equals negation. We derive exact displacement counts and, on prime-power grids with $p>5$, obtain the probability and time distribution of recovery following a random king-move perturbation. For odd primes other than five, we establish stability of orbit membership under prime-power refinement for primitive points sharing a quadratic invariant. Finally, we characterise Wall--Sun--Sun primes through directed near-return lags and, equivalently, recovery probabilities.
\end{abstract}

\maketitle
\section{Introduction}
Imagine dancers on a chequered floor. At each beat a reversible choreography assigns each dancer a square. Every dancer eventually returns to the starting square. What else could we learn by counting visits to its neighbours?

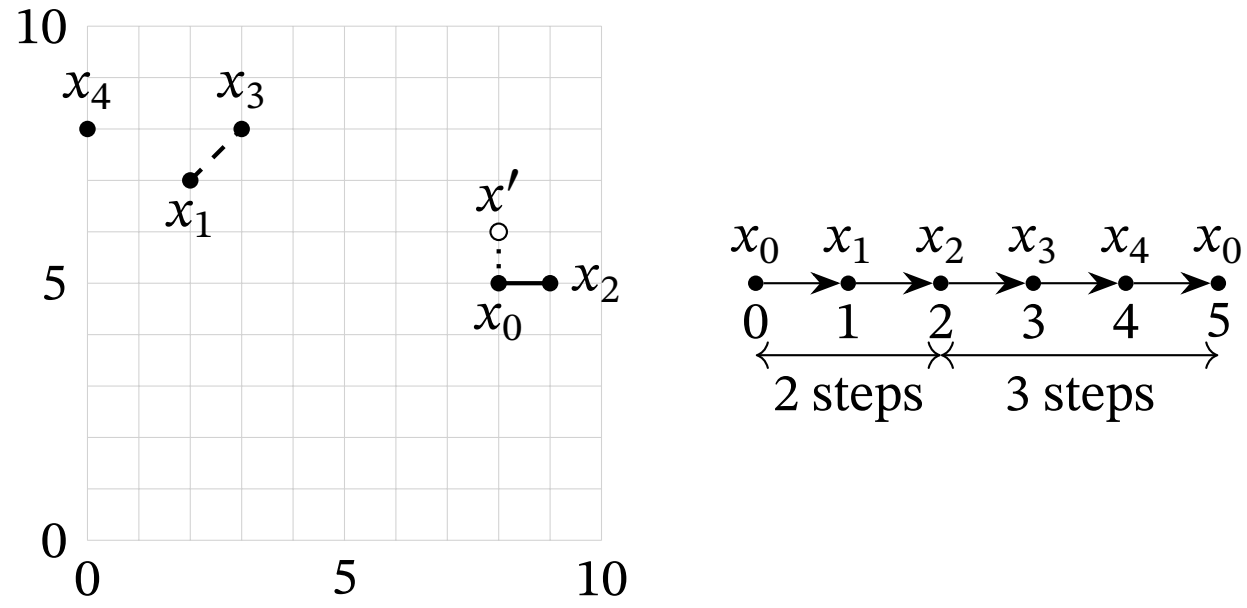
\begin{figure}[ht]
\centering
\begin{tikzpicture}[scale=.34,every node/.style={font=\small}]
\draw[step=1,gray!35,very thin] (0,0) grid (10,10);
\foreach \i in {0,5,10}{\node[below] at (\i,0){\i};\node[left] at (0,\i){\i};}
\draw[thick] (8,5)--(9,5);\draw[dashed,thick] (2,7)--(3,8);\draw[dotted,thick] (8,5)--(8,6);
\foreach \x/\y/\i/\pos in {8/5/0/below,2/7/1/below,9/5/2/right,3/8/3/above,0/8/4/above}{\fill (\x,\y) circle (1.6mm);\node[\pos] at (\x,\y){$x_\i$};}
\draw[fill=white] (8,6) circle (1.6mm);\node[above] at (8,6){$x'$};
\begin{scope}[xshift=13cm,yshift=5cm]
\foreach \i in {0,...,5}{\pgfmathtruncatemacro{\j}{mod(\i,5)}\fill (1.8*\i,0) circle (1.5mm);\node[above] at (1.8*\i,0){$x_\j$};\node[below] at (1.8*\i,0){\i};}
\foreach \i in {0,...,4}{\draw[-{Stealth[length=2mm]}] (1.8*\i+.15,0)--(1.8*\i+1.65,0);}
\draw[<->] (0,-1.4)--node[below]{$2$ steps}(3.6,-1.4);
\draw[<->] (3.6,-1.4)--node[below]{$3$ steps}(9,-1.4);
\end{scope}
\end{tikzpicture}
\caption{One cycle on side eleven. Solid and dashed segments join spatial neighbours on the cycle; the dotted segment joins points on different cycles. The time order gives two directed lags on each retained edge.}\label{fig:orbit}
\end{figure}

For $m\ge2$, let $G_m=(\Z/m\Z)^2$, a square grid with opposite sides identified, and define the discrete cat map by
\[
 A=\begin{pmatrix}1&1\\1&2\end{pmatrix},\qquad T_m(x)=Ax\pmod m.
\]
Since $\det A=1$, the inverse of $A$ has integer entries. Thus $T_m$ permutes the $m^2$ squares, and each point lies on a cycle. For example, on side eleven the orbit of $(8,5)$ is
\[
 (8,5)\to(2,7)\to(9,5)\to(3,8)\to(0,8)\to(8,5).
\]
It returns exactly after five steps and visits the neighbouring square $(9,5)$ after two steps. The neighbour $(8,6)$, however, belongs to another cycle. Figure~\ref{fig:orbit} illustrates the distinction between spatial adjacency and time order.

Cat-map periods have been studied in \cite{PV,DF,Chen}. Gaspari \cite{Gaspari} uses the quadratic invariant on prime grids. Recurrence plots and networks \cite{EOR,Donner} encode spatial proximity along trajectories. Here we retain both directed lags on each near-return edge, thereby distinguishing points with identical exact-return landscapes \cite{ALF}. Section~2 sets up the generating function, Section~3 classifies the king-move components, and Sections~4--5 count recovery times and study refinement.

\section{Near-return graphs and landscapes}
Let $T$ be a permutation of a finite set $X$ and let $k(x)$ be the length of the cycle containing $x$. A simple graph $B$ on $X$ specifies spatial adjacency.
\begin{definition}
The \emph{near-return graph} $R$ joins two points precisely when they are neighbours in $B$ and lie on the same cycle of $T$. For an edge $\{x,y\}$, the \emph{directed lag} $\ell(x,y)$ is the least positive $n$ for which $T^nx=y$.
\end{definition}
Going from $x$ to $y$ and then back to $x$ completes one circuit. Hence
\begin{equation}\label{eq:lag}
 \ell(x,y)+\ell(y,x)=k(x).
\end{equation}
Let $N_B[x]$ denote $x$ together with its neighbours in $B$, and let $N_R(x)$ denote the neighbours of $x$ in $R$ (excluding $x$). We record each visit to $N_B[x]$ at time $n$ by a term $z^n$. The arithmetic landscape is
\begin{equation}\label{eq:landscape}
 \Lscape_B(x;z)=\sum_{\substack{n\ge0\\ T^nx\in N_B[x]}}z^n,\qquad |z|<1.
\end{equation}
The sum converges by comparison with the geometric series in $|z|$. On a grid we use a symmetric observation window $V$ containing zero: the neighbourhood of $x$ is $x+V$. Write $\Lscape_V$ for this landscape. The exact window is $V=\{0\}$; the king window is $V=\{-1,0,1\}^2$.
\begin{proposition}\label{prop:land}
The visits in one circuit determine the landscape:
\begin{equation}\label{eq:H}
 H_B(x;z)=1+\sum_{y\in N_R(x)}z^{\ell(x,y)},\qquad
 \Lscape_B(x;z)=\frac{H_B(x;z)}{1-z^{k(x)}}.
\end{equation}
In particular, $H_B(x;1)=1+\deg_R(x)$ and
\[
 \lim_{z\uparrow1}(1-z)\Lscape_B(x;z)=\frac{1+\deg_R(x)}{k(x)}.
\]
\end{proposition}
\begin{proof}
In the first circuit, time zero contributes $1$ and each near-return neighbour contributes $z^{\ell(x,y)}$. Every subsequent circuit repeats these terms, shifted by $k(x)$ steps. Summing the geometric series gives \eqref{eq:H}. Setting $z=1$ in the numerator counts the visits in one circuit. For the limit, cancel $1-z$ using $1-z^k=(1-z)(1+z+\cdots+z^{k-1})$ and let $z\uparrow1$.
\end{proof}
In Figure~\ref{fig:orbit}, the king landscapes at $x_0$ and $x_2$ are respectively $(1+z^2)/(1-z^5)$ and $(1+z^3)/(1-z^5)$. The exact landscape $1/(1-z^{k(x)})$ is the same for both. The limit in Proposition~\ref{prop:land} is the fraction of time spent in the starting neighbourhood.

\smallskip\noindent\textbf{Recovery after a misstep.}
Move a point from $x$ to a uniformly chosen neighbour $y$ in $B$, then iterate $T$ from $y$. We call the misstep recoverable if the orbit visits $x$ again, and write $D_x$ for the first positive time at which this happens ($D_x=\infty$ if it never happens). Recovery is possible exactly for the neighbours in $R$. By \eqref{eq:lag}, the recovery time for such a neighbour is $k(x)-\ell(x,y)$. Consequently, whenever $\deg_B(x)>0$,
\begin{equation}\label{eq:local}
 \Prob(D_x<\infty)=\frac{\deg_R(x)}{\deg_B(x)}.
\end{equation}
For example, only one of the eight king moves from $x_0$ in Figure~\ref{fig:orbit} is recoverable, and its recovery time is three. Recovery here means revisiting the original square, rather than synchronising with the unperturbed trajectory.

\section{The chess-king near-return graph}
Write $\R_m$ for the near-return graph of $T_m$ with the king window. For $m\ge7$ each square has eight king neighbours. Our first main result classifies all components of $\R_m$.
\begin{theorem}\label{thm:shapes}
Let $m\ge7$.
\begin{enumerate}
\item Every component of $\R_m$ is an isolated vertex, a single edge, a two-edge path or a triangle with a leaf.
\item There are exactly four triangles with a leaf if $5\nmid m$ and $A^n\equiv-I\pmod m$ for some positive integer $n$, and none otherwise. Their degree-three vertices are the quarter turns of $c^*=(1,-2)/5$, with division modulo $m$.
\item If $m$ is odd and $5\nmid m$, there are exactly four two-edge paths.
\end{enumerate}
\end{theorem}
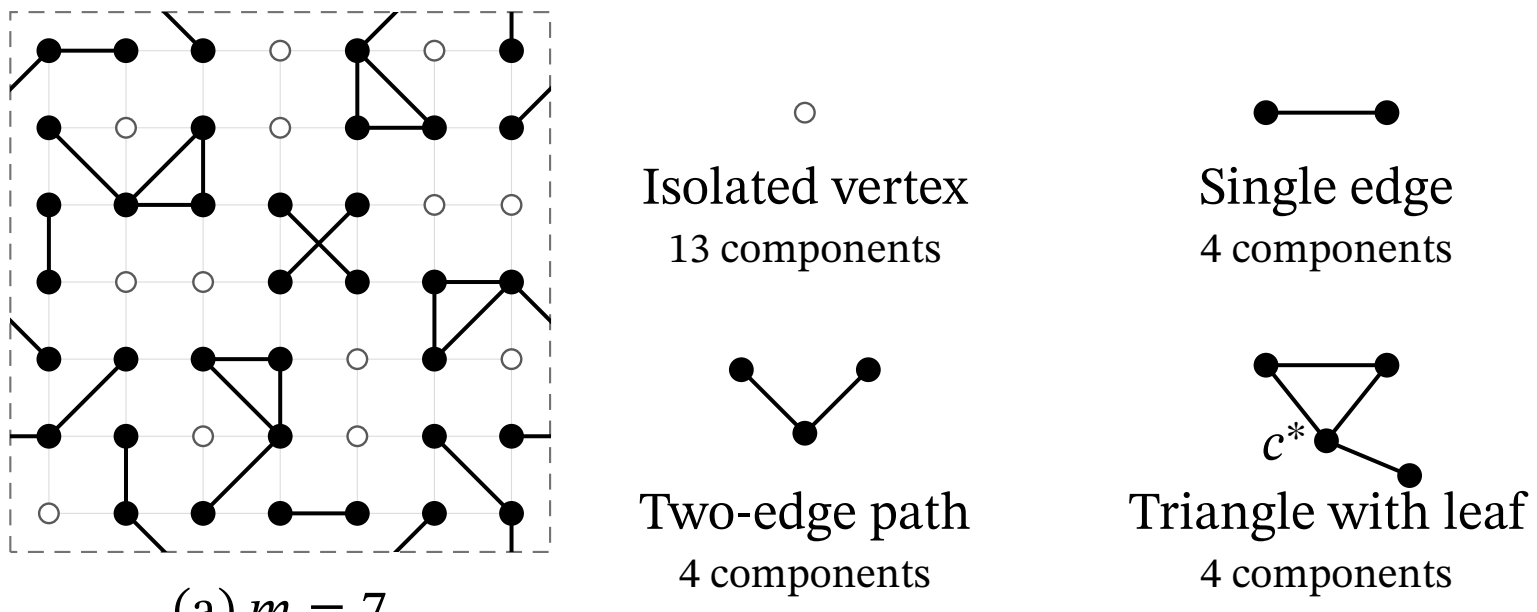
\begin{figure}[ht]
\centering
\begin{tikzpicture}[font=\small,edge/.style={line width=.75pt,line cap=round},
  dot/.style={circle,fill=black,inner sep=1.65pt,outer sep=0pt},
  isolated/.style={circle,draw=black!65,fill=white,line width=.45pt,inner sep=1.3pt,outer sep=0pt}]
\begin{scope}[x=.51cm,y=.51cm]
\draw[step=1,black!14,very thin] (0,0) grid (6,6);
\draw[densely dashed,black!55,thin] (-.5,-.5) rectangle (6.5,6.5);
\begin{scope}
\clip (-.5,-.5) rectangle (6.5,6.5);
\draw[edge] (0,1)--(1,2);
\draw[edge] (0,3)--(0,4);
\draw[edge] (0,5)--(1,4);
\draw[edge] (0,6)--(1,6);
\draw[edge] (1,0)--(1,1);
\draw[edge] (1,0)--(2,-1);
\draw[edge] (1,7)--(2,6);
\draw[edge] (1,4)--(2,4);
\draw[edge] (1,4)--(2,5);
\draw[edge] (2,0)--(3,1);
\draw[edge] (2,2)--(3,2);
\draw[edge] (2,2)--(3,1);
\draw[edge] (2,4)--(2,5);
\draw[edge] (3,0)--(4,0);
\draw[edge] (3,1)--(3,2);
\draw[edge] (3,3)--(4,4);
\draw[edge] (3,4)--(4,3);
\draw[edge] (4,5)--(5,5);
\draw[edge] (4,5)--(4,6);
\draw[edge] (4,-1)--(5,0);
\draw[edge] (4,6)--(5,7);
\draw[edge] (4,6)--(5,5);
\draw[edge] (5,1)--(6,0);
\draw[edge] (5,2)--(5,3);
\draw[edge] (5,2)--(6,3);
\draw[edge] (5,3)--(6,3);
\draw[edge] (-1,1)--(0,1);
\draw[edge] (6,1)--(7,1);
\draw[edge] (-1,3)--(0,2);
\draw[edge] (6,3)--(7,2);
\draw[edge] (-1,5)--(0,6);
\draw[edge] (6,5)--(7,6);
\draw[edge] (6,-1)--(6,0);
\draw[edge] (6,6)--(6,7);
\end{scope}
\node[isolated] at (0,0) {};
\node[dot] at (0,1) {};
\node[dot] at (0,2) {};
\node[dot] at (0,3) {};
\node[dot] at (0,4) {};
\node[dot] at (0,5) {};
\node[dot] at (0,6) {};
\node[dot] at (1,0) {};
\node[dot] at (1,1) {};
\node[dot] at (1,2) {};
\node[isolated] at (1,3) {};
\node[dot] at (1,4) {};
\node[isolated] at (1,5) {};
\node[dot] at (1,6) {};
\node[dot] at (2,0) {};
\node[isolated] at (2,1) {};
\node[dot] at (2,2) {};
\node[isolated] at (2,3) {};
\node[dot] at (2,4) {};
\node[dot] at (2,5) {};
\node[dot] at (2,6) {};
\node[dot] at (3,0) {};
\node[dot] at (3,1) {};
\node[dot] at (3,2) {};
\node[dot] at (3,3) {};
\node[dot] at (3,4) {};
\node[isolated] at (3,5) {};
\node[isolated] at (3,6) {};
\node[dot] at (4,0) {};
\node[isolated] at (4,1) {};
\node[isolated] at (4,2) {};
\node[dot] at (4,3) {};
\node[dot] at (4,4) {};
\node[dot] at (4,5) {};
\node[dot] at (4,6) {};
\node[dot] at (5,0) {};
\node[dot] at (5,1) {};
\node[dot] at (5,2) {};
\node[dot] at (5,3) {};
\node[isolated] at (5,4) {};
\node[dot] at (5,5) {};
\node[isolated] at (5,6) {};
\node[dot] at (6,0) {};
\node[dot] at (6,1) {};
\node[isolated] at (6,2) {};
\node[dot] at (6,3) {};
\node[isolated] at (6,4) {};
\node[dot] at (6,5) {};
\node[dot] at (6,6) {};
\node[below=3pt] at (3,-.5) {(a) $m=7$};
\end{scope}
\begin{scope}[shift={(5.0cm,0cm)}]
\node[isolated] at (0,2.65) {};
\node[align=center] at (0,1.95) {Isolated vertex\\[-1pt]\scriptsize 13 components};
\draw[edge] (3.05,2.65)--(3.85,2.65);
\foreach \p in {(3.05,2.65),(3.85,2.65)}{\node[dot] at \p {};}
\node[align=center] at (3.45,1.95) {Single edge\\[-1pt]\scriptsize 4 components};
\draw[edge] (-.42,.95)--(0,.53)--(.42,.95);
\foreach \p in {(-.42,.95),(0,.53),(.42,.95)}{\node[dot] at \p {};}
\node[align=center] at (0,-.20) {Two-edge path\\[-1pt]\scriptsize 4 components};
\draw[edge] (3.05,.98)--(3.85,.98)--(3.45,.48)--cycle;
\draw[edge] (3.45,.48)--(4.0,.25);
\foreach \p in {(3.05,.98),(3.85,.98),(3.45,.48),(4.0,.25)}{\node[dot] at \p {};}
\node[left,inner sep=3pt] at (3.45,.48) {$c^*$};
\node[align=center] at (3.45,-.20) {Triangle with leaf\\[-1pt]\scriptsize 4 components};
\node at (1.725,-.85) {(b) Component types};
\end{scope}
\end{tikzpicture}
\caption{The king near-return graph on the torus of side seven and its four component types. Edges meeting the dashed boundary continue across the opposite boundary. Open vertices are isolated. The component counts give 49 vertices and 28 edges. In the last schematic, $c^*$ denotes a degree-three vertex up to a quarter turn.}\label{fig:shapes}
\end{figure}
\subsection{The quadratic invariant}
For $x=(x_1,x_2)$ put
\[
 Q(x)=x_2^2-x_1x_2-x_1^2=-\det(x,Ax).
\]
Since $\det A=1$, we have $\det(Ax,A^2x)=\det(x,Ax)$, so $Q(Ax)=Q(x)$. In particular, points on the same orbit have equal $Q$ modulo $m$. Call a king edge a \emph{candidate edge} if its endpoints have equal $Q$. Every near-return edge is a candidate edge; the converse need not hold.

Put $U=2x_2-x_1$ and $W=x_2+2x_1$. For a displacement $v=(v_1,v_2)$, expansion gives
\begin{equation}\label{eq:diff}
 Q(x+v)-Q(x)=v_2U-v_1W+Q(v).
\end{equation}
For instance, the upward move $(0,1)$ preserves $Q$ exactly when $U=-1$. Applying \eqref{eq:diff} to all eight moves gives Figure~\ref{fig:tests}, or, grouping opposite directions,
\begin{equation}\label{eq:tests}
 U=\pm1,\qquad W=\pm1,\qquad U-W=\pm1,\qquad U+W=\pm1.
\end{equation}
All congruences in this section are modulo $m$. Each of the four families in \eqref{eq:tests} contributes at most one candidate edge.
\begin{figure}[ht]
\centering
\begin{tikzpicture}[x=2.65cm,y=.58cm]
\foreach \i in {-1.5,-.5,.5,1.5}{\draw (\i,-1.5)--(\i,1.5);\draw (-1.5,\i)--(1.5,\i);}
\foreach \x/\y/\t in {-1/1/{U+W=-1},0/1/{U=-1},1/1/{U-W=1},-1/0/{W=1},0/0/{x},1/0/{W=-1},-1/-1/{U-W=-1},0/-1/{U=1},1/-1/{U+W=1}}{\node at (\x,\y){$\t$};}
\end{tikzpicture}
\caption{Candidate-edge tests for the eight king moves. Columns and rows correspond to displacements $-1,0,1$ in the first and second coordinates, respectively.}\label{fig:tests}
\end{figure}

The quarter turn $J(x_1,x_2)=(-x_2,x_1)$ preserves king adjacency. Substitution gives $Q(Jx)=-Q(x)$ and $JAJ^{-1}=A^{-1}$. Thus rotation preserves both candidate and near-return edges, reversing the time order of each orbit. In the new coordinates it sends $(U,W)$ to $(W,-U)$.

\begin{lemma}\label{lem:centres}
If a vertex has at least two candidate edges, its pair $(U,W)$ is a quarter turn of one of
\begin{equation}\label{eq:centres}
 (1,1),\quad(-1,0),\quad(2,-1),\quad(1,-2).
\end{equation}
For even $m$ there are also the quarter turns of $(m/2-1,m/2)$.
\end{lemma}
\begin{proof}
Choose two candidate moves. They belong to distinct families in \eqref{eq:tests}.

\noindent\textit{Case 1. Two axial moves.} Then $U=\pm1$ and $W=\pm1$. These are exactly the quarter turns of $(1,1)$.

\noindent\textit{Case 2. One axial and one diagonal move.} By rotation we may assume $U=1$. The diagonal test is $1-W=\pm1$ or $1+W=\pm1$, giving $W=0,2,-2$. These pairs are turns of the last three representatives in \eqref{eq:centres}.

\noindent\textit{Case 3. Two diagonal moves.} Write $\alpha=U-W$ and $\beta=U+W$. A quarter turn sends $(\alpha,\beta)$ to $(\beta,-\alpha)$, so we may take $\alpha=\beta=-1$. The equations give $U=W-1$ and $2W=0$. If $m$ is odd, $W=0$; if $m$ is even, $W=0$ or $m/2$. These give the two stated representatives.
\end{proof}

\subsection{Determining which edges survive}
\begin{lemma}\label{lem:triangle}
Suppose $5\nmid m$. At each quarter turn of $c^*=(1,-2)/5$, the candidate component is a triangle with a leaf. These four components are distinct. Each retains its shape in $\R_m$ if a power of $A$ equals $-I$ modulo $m$, and otherwise splits into two single edges.
\end{lemma}
\begin{proof}
The change $x\mapsto(U,W)$ has determinant $-5$ and is therefore invertible. The pair $(-1,0)$ corresponds to $c^*$. Applying \eqref{eq:tests} at this point and its neighbours gives the following table:
\[
\begin{array}{c|c|c|c}
\text{Move }v&c^*+v&\text{New }(U,W)&\text{Degree}\\\hline
(0,1)&-Ac^*&(1,1)&2\\
(-1,1)&-A^{-1}c^*&(2,-1)&2\\
(-1,-1)&A^2c^*&(-2,-3)&1
\end{array}
\]
The first two endpoints are neighbours with equal $Q$; they form a triangle with $c^*$. The third endpoint is its leaf. The degrees in the table exclude any further edges. Rotation gives four distinct components since their unique degree-three vertices are distinct.

The edges $\{c^*,A^2c^*\}$ and $\{-A^{-1}c^*,-Ac^*\}$ always survive. The other two join the orbit of $c^*$ to that of $-c^*$, so both survive exactly when some power $A^n$ sends $c^*$ to $-c^*$.

\emph{Sufficiency.} If $A^n=-I$, then $A^nc^*=-c^*$.

\emph{Necessity.} Suppose $A^nc^*=-c^*$. Then $A^n(Ac^*)=-Ac^*$ as well. Since $\det(c^*,Ac^*)=-1/5$ is invertible modulo $m$, every vector is a linear combination of $c^*$ and $Ac^*$. The matrix $A^n$ negates both, hence $A^n=-I$. If no such power exists, the two other edges are deleted and the component splits into two edges.
\end{proof}
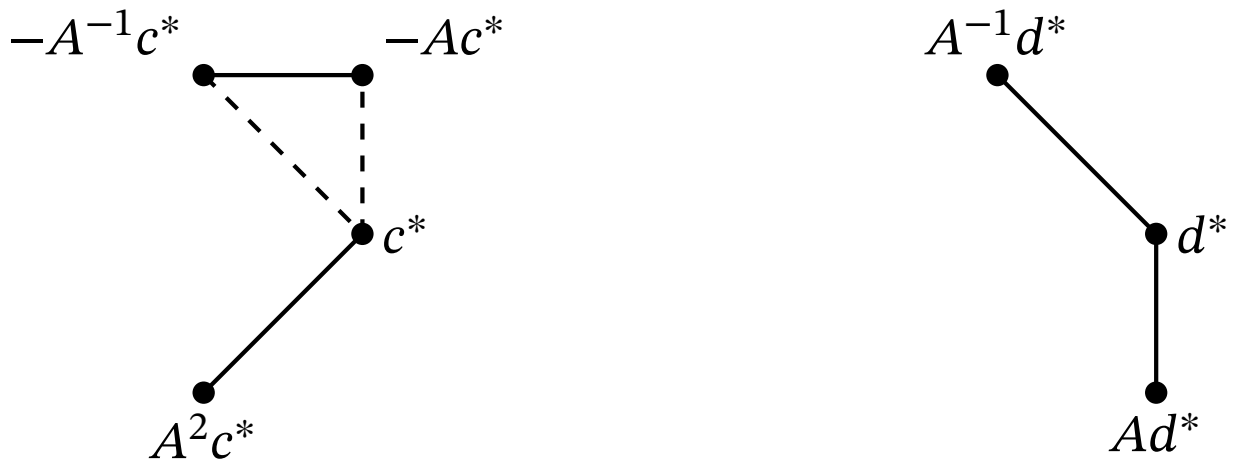
\begin{figure}[ht]
\centering
\begin{tikzpicture}[scale=1.05,every node/.style={font=\small}]
\draw[thick] (0,1)--(1,1) (1,0)--(0,-1);\draw[thick,dashed] (0,1)--(1,0)--(1,1);
\foreach \p in {(0,1),(1,1),(1,0),(0,-1)}{\fill \p circle(2pt);}
\node[above left] at (0,1){$-A^{-1}c^*$};\node[above right] at (1,1){$-Ac^*$};\node[right] at (1,0){$c^*$};\node[below] at (0,-1){$A^2c^*$};
\begin{scope}[xshift=5cm]
\draw[thick] (0,1)--(1,0)--(1,-1);\foreach \p in {(0,1),(1,0),(1,-1)}{\fill \p circle(2pt);}
\node[above] at (0,1){$A^{-1}d^*$};\node[right] at (1,0){$d^*$};\node[below] at (1,-1){$Ad^*$};
\end{scope}
\end{tikzpicture}
\caption{The candidate components at $c^*$ and $d^*=(-1,0)$. Solid edges always survive; dashed edges survive exactly when a power of $A$ equals $-I$. The path at $d^*$ is a component when $5\nmid m$, including $m=8$ after the extra candidate edge is deleted.}\label{fig:components}
\end{figure}

\begin{proof}[Proof of Theorem~\ref{thm:shapes}]
Any component containing at least three vertices has a vertex of degree at least two. It therefore suffices to examine the families in Lemma~\ref{lem:centres}.

First, at the representative $(1,-2)$, the two candidate moves and their endpoint coordinates are
\[
\begin{array}{c|c|c}
 v&\text{New }(U,W)&\text{Degree if }m\ne8\\\hline
 (0,-1)&(-1,-3)&1\\
 (-1,1)&(4,-3)&1
\end{array}
\]
The degrees follow from \eqref{eq:tests}. Thus this representative gives a two-edge candidate path when $m\ne8$.

\noindent\textit{Case 1. $5\mid m$.} Every grid point satisfies $2U+W=5x_2\equiv0\pmod5$. To test all rotations of each representative, it suffices to check $2U+W$ and $2W-U$, since the other two rotations negate these values:
\[
\begin{array}{c|rrrr}
 (U,W)&(1,1)&(-1,0)&(2,-1)&(1,-2)\\\hline
 2U+W&3&-2&3&0\\
 2W-U&1&1&-4&-5
\end{array}
\]
Only the last family is possible. If $m$ is even, the extra representative reduces to $(-1,0)$ modulo five and is excluded by the same test. Therefore every candidate component with at least three vertices is a two-edge path. Deleting edges can only shorten or disconnect such a path.

\noindent\textit{Case 2. $5\nmid m$, $m\ne8$.} The change of coordinates is invertible. Lemma~\ref{lem:triangle} accounts for the first three families in \eqref{eq:centres}. The fourth corresponds to $d^*=(-1,0)$, whose endpoints are $Ad^*=(-1,-1)$ and $A^{-1}d^*=(-2,1)$. Both edges survive, and rotation gives four distinct two-edge paths (Figure~\ref{fig:components}).

For odd $m$ this exhausts the possibilities. For even $m$, write $h=m/2>4$. The extra centre $(h-1,h)$ gives two candidate moves:
\[
\begin{array}{c|c|c}
 v&\text{New }(U,W)&\text{Test passed}\\\hline
 (-1,1)&(h+2,h-1)&U+W=1\\
 (-1,-1)&(h-2,h-3)&U-W=1
\end{array}
\]
At each endpoint neither coordinate is $\pm1$. The other diagonal values are $3$ and $-5$, respectively, neither of which is $\pm1$ modulo $m\ge10$. Thus both endpoints have degree one. These extra components are candidate paths; the orbit test can only delete edges.

\noindent\textit{Case 3. $m=8$.} Here $A$ has order six and $A^3=5I\not\equiv-I\pmod8$, so the four candidate triangles split into edges. In the remaining families the additional centres are the turns of $(4,-3)$, one of which corresponds to $A^{-1}d^*=(6,1)$. The orbit of $d^*=(7,0)$ is
\[
 (7,0)\to(7,7)\to(6,5)\to(3,0)\to(3,3)\to(6,1)\to(7,0).
\]
The only king neighbour of $(6,1)$ on this orbit is $(7,0)$. Hence the extra candidate edge at $(6,1)$ is deleted, and so are its rotations. The two edges incident to $d^*$ survive, with the other endpoint of candidate degree one by \eqref{eq:tests}. Thus these components are two-edge paths.

Every remaining component has maximum degree one and is therefore a vertex or an edge. Lemma~\ref{lem:triangle} gives the triangle count, and Case 2 gives exactly four two-edge paths when $m$ is odd and $5\nmid m$.
\end{proof}

The bound $m\ge7$ is sharp: labelling the cycles and retaining king edges between equally labelled points gives components of ten vertices on side five and six vertices on side six.

\subsection{Counting edges and recoverable moves}
Theorem~\ref{thm:shapes} and \eqref{eq:local} imply that the probability of recovery from a uniformly chosen king move at a fixed point is one of $0,1/8,1/4,3/8$. The value $3/8$ occurs at exactly four points when the triangle condition holds, and nowhere otherwise.

There are exactly $4m$ candidate edges. Indeed, choose the four directions $(1,0)$, $(0,1)$, $(1,1)$, $(1,-1)$, counting each unoriented edge once. For each direction \eqref{eq:diff} is a linear congruence in $x_1,x_2$ with one coefficient equal to $\pm1$. Choosing the other coordinate freely gives exactly $m$ solutions. Therefore $|E(\R_m)|\le4m$. If both the start and move are chosen independently and uniformly, there are $8m^2$ outcomes, and each edge gives two successful outcomes. Hence
\begin{equation}\label{eq:global}
 \Prob(D<\infty)=\frac{2|E(\R_m)|}{8m^2}\le\frac1m.
\end{equation}
Equality holds exactly when every candidate edge survives. For example, all 28 candidates survive on side seven, but only 16 of 44 survive on side eleven. The next section determines this edge count on prime-power grids.

\section{Counting visits and recovery times}
Fix an odd prime $p\ne5$. We first determine when a prescribed displacement $v$ can occur after $n$ iterations on a grid of side $p^a$. Call a vector \emph{primitive modulo $p$} if at least one coordinate is nonzero modulo $p$. Every nonzero king displacement is primitive.

\subsection{Fibonacci arithmetic}
Let $F_0=0$, $F_1=1$ and $F_{j+2}=F_{j+1}+F_j$. Put
\[
 M=\begin{pmatrix}0&1\\1&1\end{pmatrix},\qquad H=2M-I.
\]
Direct multiplication gives $M^2=A$, $M^2=M+I$, $H^2=5I$ and $A=(3I+H)/2$. Induction using $M^2=M+I$ gives
\begin{equation}\label{eq:fib}
 M^j=F_jM+F_{j-1}I=\begin{pmatrix}F_{j-1}&F_j\\F_j&F_{j+1}\end{pmatrix}\quad(j\ge1).
\end{equation}
Let $\varepsilon_p=(5/p)$ be the Legendre symbol: it is $1$ if five is a square modulo $p$, and $-1$ otherwise. Euler's criterion gives $5^{(p-1)/2}\equiv\varepsilon_p\pmod p$. Define
\[
 \tau_p=\ord_p(A),\qquad s_p=\nu_p(F_{p-\varepsilon_p}),
\]
where $\ord_p(A)$ is the least positive $n$ for which $A^n\equiv I\pmod p$, and $\nu_p(b)$ is the exponent of $p$ in $b$ (for nonzero rationals, the exponent in the numerator minus that in the denominator). Write $K_p(a)=\ord_{p^a}(A)$. A matrix with determinant prime to $p$ is invertible modulo every $p^a$, by the two-by-two inverse formula.

\begin{lemma}\label{lem:arithmetic}
Let $p$ be an odd prime other than five. Then $\tau_p\mid p-\varepsilon_p$ and $s_p\ge1$. For positive $n$:
\begin{enumerate}
\item If $\tau_p\nmid n$, then $A^n-I$ is invertible modulo $p$.
\item If $n=\tau_p j$, then $A^n-I=p^{s_p+\nu_p(j)}B_j$ for an integer matrix $B_j$ invertible modulo $p$.
\end{enumerate}
Consequently,
\begin{equation}\label{eq:order}
 K_p(a)=\tau_p p^{\max(0,a-s_p)},
\end{equation}
and every primitive point modulo $p$ has this period on side $p^a$.
\end{lemma}
\begin{proof}
\emph{The order modulo $p$.} The binomial theorem and Fermat's little theorem give
\[
 A^p\equiv\frac{3I+H^p}{2}
 =\frac{3I+5^{(p-1)/2}H}{2}
 \equiv\frac{3I+\varepsilon_pH}{2}\pmod p.
\]
If $\varepsilon_p=1$, this equals $A$; otherwise it equals $A^{-1}$, since $(3I+H)(3I-H)=4I$. In both cases $A^{p-\varepsilon_p}\equiv I$, so $\tau_p\mid p-\varepsilon_p$.

\emph{Part 1.} We prove the contrapositive. Since $H^2=5I$, every power of $A$ is of the form $cI+dH$ modulo $p$. Its determinant is $c^2-5d^2$. If $A^n-I$ is singular modulo $p$, then
\[
 c^2-5d^2=1,\qquad(c-1)^2-5d^2=0.
\]
Subtracting gives $c=1$ and hence $5d^2=0$. As $p\ne2,5$, we get $d=0$. Thus $A^n=I$ modulo $p$, which implies $\tau_p\mid n$.

\emph{Part 2.} First determine the common power of $p$ in the entries of $A^{\tau_p}-I$. Put $h=p-\varepsilon_p$. Since $h$ is even, $\det M^h=1$. The two-by-two inverse formula gives $M^{-h}=\operatorname{tr}(M^h)I-M^h$. Using \eqref{eq:fib}, we obtain $M^h-M^{-h}=F_hH$. Multiplying by $M^h$ yields
\begin{equation}\label{eq:Fidentity}
 A^h-I=F_h M^hH.
\end{equation}
The determinant of $M^hH$ is $-5$, prime to $p$. Hence $s_p\ge1$, and after division by $p^{s_p}$ the matrix in \eqref{eq:Fidentity} is invertible modulo $p$.

Write $h=q\tau_p$. Factoring,
\[
 A^h-I=(A^{\tau_p}-I)S,\qquad S=I+A^{\tau_p}+\cdots+A^{(q-1)\tau_p}.
\]
Here $S\equiv qI\pmod p$ and $p\nmid q$ since $q\mid p-\varepsilon_p$. Cancelling the invertible matrix $S$ modulo $p^{s_p}$ shows $A^{\tau_p}=I+p^{s_p}B$. Dividing the factorisation by $p^{s_p}$ and using \eqref{eq:Fidentity} shows that $B$ is invertible modulo $p$.

Now write $j=p^ru$ with $p\nmid u$. For $t\ge1$, the binomial theorem gives
\begin{equation}\label{eq:lte}
 (I+p^tC)^p\equiv I+p^{t+1}C\pmod{p^{t+2}}.
\end{equation}
Indeed, each middle term contains $p^{2t+1}$, and the last term contains $p^{tp}$; both exponents are at least $t+2$ for odd $p$. Repeating \eqref{eq:lte} $r$ times gives $(I+p^{s_p}B)^{p^r}=I+p^{s_p+r}C_r$ with $C_r\equiv B\pmod p$. Raising to $u$, and expanding modulo $p^{s_p+r+1}$, leaves the linear term $up^{s_p+r}C_r$. Thus
\[
 A^{\tau_p j}-I=p^{s_p+r}B_j,\qquad B_j\equiv uB\pmod p.
\]
The last matrix is invertible modulo $p$, as required.

\emph{Periods.} By Parts 1 and 2, $A^n\equiv I\pmod{p^a}$ holds exactly when $n=\tau_p j$ and $s_p+\nu_p(j)\ge a$. The least such $j$ is $p^{\max(0,a-s_p)}$, proving \eqref{eq:order}. For a primitive $x$, the same condition is necessary and sufficient for $(A^n-I)x=0$ modulo $p^a$: in Part 1 invertibility excludes a nonzero reduction, and in Part 2 the vector $B_jx$ is still primitive. Thus the period of $x$ is $K_p(a)$.
\end{proof}
The period formula \eqref{eq:order} is the familiar Fibonacci period-lifting law \cite{Wall,DF,Chen}; we have included the matrix argument to make explicit the additional invertibility used here. The period is constant through level $s_p$, and then grows by a factor of $p$ at each refinement.

\subsection{Counting displacements}
For a positive lag $n$ and displacement $v$, put
\[
 \Vis_m(n,v)=\#\{x\in G_m:A^nx=x+v\}.
\]
The visit equation is $(A^n-I)x=v$. For an odd prime $p\ne5$, set
\[
 d=\begin{cases}0,&\tau_p\nmid n,\\ \min\{a,s_p+\nu_p(n/\tau_p)\},&\tau_p\mid n.\end{cases}
\]
\begin{equation}\label{eq:allvis}
 \Vis_{p^a}(n,v)=\begin{cases}p^{2d},&p^d\mid v_1\text{ and }p^d\mid v_2,\\0,&\text{otherwise}.\end{cases}
\end{equation}
To prove \eqref{eq:allvis}, consider two cases.

\noindent\textit{Case 1. $\tau_p\nmid n$.} Lemma~\ref{lem:arithmetic} makes $A^n-I$ invertible, giving exactly one solution, as required for $d=0$.

\noindent\textit{Case 2. $\tau_p\mid n$.} Put $t=s_p+\nu_p(n/\tau_p)$. The same lemma gives $A^n-I=p^t B$ with $B$ invertible. Substituting $u=Bx$ reduces the equation to $p^t u_i\equiv v_i\pmod{p^a}$ for $i=1,2$. If $t<a$, a solution exists precisely when $p^t\mid v_i$; dividing by $p^t$ fixes $u_i$ modulo $p^{a-t}$ and gives $p^t$ choices modulo $p^a$. If $t\ge a$, every residue works when $v_i=0$ modulo $p^a$, and none otherwise. Thus each coordinate has $p^d$ choices when $p^d\mid v_i$, and none otherwise. Multiplying the two counts proves \eqref{eq:allvis}. For example, on side $49$ at lag $8$, we have $\tau_7=8$, $s_7=1$ and $d=1$: there are $49$ exact returns, $49$ visits with displacement $(7,0)$, and none with displacement $(1,0)$.
\begin{theorem}\label{thm:recovery}
Let $p>5$ be prime, $a\ge1$, $K=K_p(a)$ and $\tau=\tau_p$. For every primitive displacement $v$ on $G_{p^a}$,
\begin{equation}\label{eq:visit}
 \Vis_{p^a}(n,v)=\begin{cases}1,&\tau\nmid n,\\0,&\tau\mid n.\end{cases}
\end{equation}
Choose the starting point and a nonzero king move independently and uniformly. The first recovery time $D$ satisfies
\begin{equation}\label{eq:dist}
 \Prob(D=n)=\begin{cases}p^{-2a},&1\le n<K,\ \tau\nmid n,\\0,&\text{otherwise},\end{cases}
\end{equation}
for positive $n$. In particular,
\begin{equation}\label{eq:recover}
 \begin{aligned}
 |E(\R_{p^a})|&=4K(1-1/\tau),\\
 \Prob(D<\infty)&=\frac{K(1-1/\tau)}{p^{2a}},&
 \E[D\mid D<\infty]&=\frac K2.
 \end{aligned}
\end{equation}
\end{theorem}
\begin{proof}
For a primitive displacement, $d=0$ is the only case in which the divisibility condition in \eqref{eq:allvis} holds. This gives \eqref{eq:visit}.

For recovery, let $x$ be the starting point and $y=x+v$ the displaced point. Returning to $x$ at time $n$ requires
\[
 A^ny=y-v,\qquad\text{or}\qquad (A^n-I)y=-v.
\]
Every solution $y$ is primitive, since otherwise the left-hand side would vanish modulo $p$. Its period is therefore $K$ by Lemma~\ref{lem:arithmetic}. As $x\ne y$, the first recovery occurs at a unique lag in $1\le n<K$. For each such lag with $\tau\nmid n$, \eqref{eq:visit} gives one starting pair per direction. This is eight successful outcomes out of $8p^{2a}$ equally likely choices, proving \eqref{eq:dist}.

There are $K-K/\tau$ admissible lags. Summing \eqref{eq:dist} gives the recovery probability; each undirected edge accounts for two recoverable starting pairs, giving the edge count. Conditional on recovery all admissible lags are equally likely. Pairing each lag $n$ with $K-n$ preserves the admissible set, so their mean is $K/2$.
\end{proof}
The cumulative probability for $0\le h<K$ is $(h-\lfloor h/\tau\rfloor)/p^{2a}$, by counting the admissible lags up to $h$. The argument works for any fixed primitive displacement; the distribution therefore also holds for a window of such moves chosen independently of the starting point.
\begin{table}[ht]
\centering
\caption{Recovery statistics on selected prime-power grids.}\label{tab:recovery}
\begin{tabular}{rrrrr}\toprule
Side&$K$&Edges&$\Prob(D<\infty)$&$\E[D\mid D<\infty]$\\\midrule
7&8&28&$1/7$&4\\
11&5&16&$4/121$&$5/2$\\
13&14&52&$1/13$&7\\
49&56&196&$1/49$&28\\
121&55&176&$4/1331$&$55/2$\\\bottomrule
\end{tabular}
\end{table}
Table~\ref{tab:recovery} illustrates the difference between local and global recovery. On side seven the global bound in \eqref{eq:global} is attained, while on side eleven it is not.

\section{Lifting near returns}
On refining from side $p^a$ to side $p^{a+1}$, an orbit may make several circuits of the coarse orbit before closing. We ask whether a given pair of endpoints stays on a common orbit, and which circuit reaches the chosen fine endpoint.

\begin{figure}[ht]
\centering
\begin{tikzpicture}[every node/.style={font=\small}]
\begin{scope}[xshift=0cm]
\draw (0,0) circle (.88cm);
\foreach \i in {0,...,3}{\pgfmathsetmacro{\ang}{90-90*\i}\fill (\ang:.88) circle(2pt);\node at (\ang:1.15){\i};}
\draw[thick] (180:.88) circle(4pt);\node[below] at (0,-1.45){Side three};
\end{scope}
\draw[-{Stealth}] (1.45,0)--(2.4,0);
\begin{scope}[xshift=4.1cm]
\draw (0,0) circle (1.2cm);
\foreach \i in {0,...,11}{\pgfmathsetmacro{\ang}{90-30*\i}\fill (\ang:1.2) circle(2pt);\node at (\ang:1.46){\i};}
\foreach \ang in {0,-120,-240}{\draw[thick] (\ang:1.2) circle(4pt);}
\draw[thick] (-120:1.2) circle(6pt);\node[below] at (0,-1.65){Side nine};
\end{scope}
\node[align=left,anchor=west] at (5.95,0){$3,7,11$ reduce\\to step $3$;\\step $7$ hits\\the target.};
\end{tikzpicture}
\caption{One extra circuit of the coarse orbit is required to reach the chosen fine neighbour.}\label{fig:lift}
\end{figure}
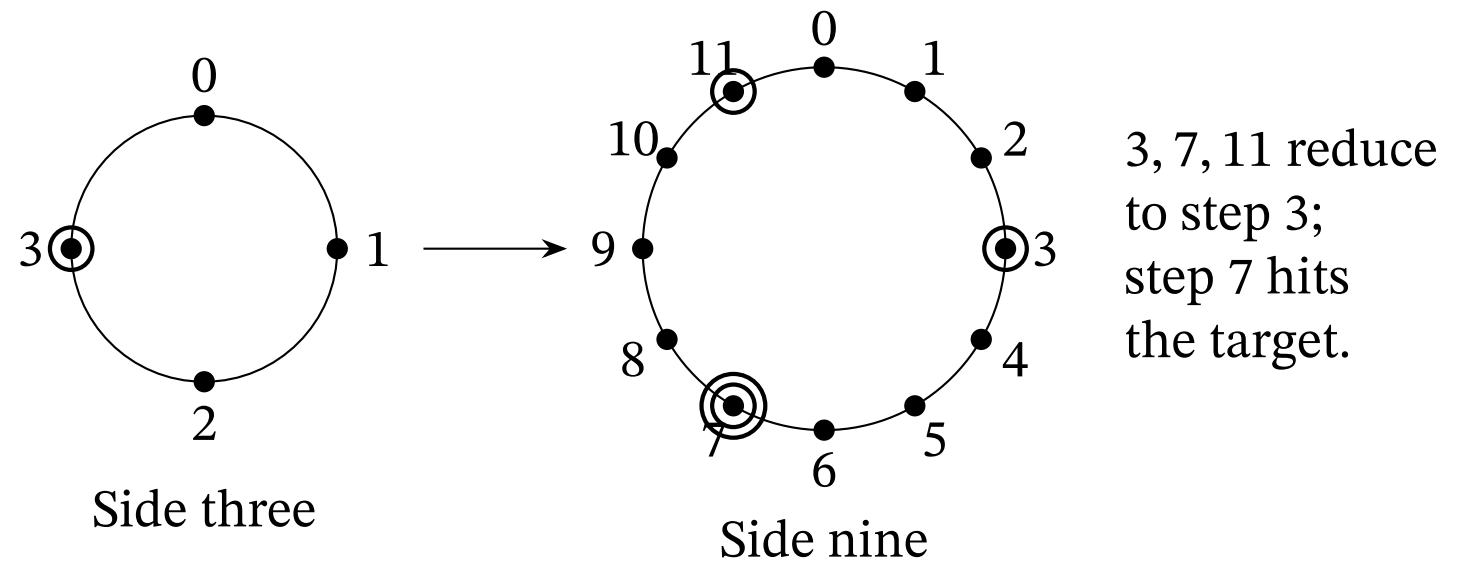

For example, the edge $c^*\to c^*+(0,1)$ from Section~3 lifts from side three to side nine. The coarse lag is three and the period is four. The fine period is twelve. The points reached at steps $3,7,11$ all reduce to the same coarse endpoint, but only step seven reaches the chosen neighbour (Figure~\ref{fig:lift}).

\subsection{When does a near return persist?}
Let $p$ be an odd prime other than five. Take rational vectors $x,y$ whose coordinate denominators are prime to $p$. Their reductions modulo $p^a$, obtained by inverting these denominators, are denoted $x_a,y_a$.
\begin{theorem}\label{thm:lift}
Suppose $x$ is primitive modulo $p$ and $Q(y)=Q(x)$ over the rationals. For every $a\ge s_p$, the points $x_a,y_a$ lie on the same orbit if and only if $x_{a+1},y_{a+1}$ do. If $0\le n<K_p(a)$ is the first visit time from $x_a$ to $y_a$ (allowing zero), the fine visit time is
\[
 n+jK_p(a)\quad\text{for a unique }j\in\{0,1,\ldots,p-1\}.
\]
\end{theorem}
\begin{proof}
\emph{Necessity.} If $x_{a+1},y_{a+1}$ lie on the same orbit, reducing a corresponding iterate modulo $p^a$ gives a visit from $x_a$ to $y_a$.

\emph{Sufficiency.} Suppose $A^nx_a=y_a$ at time $0\le n<K_p(a)$. Write $K=K_p(a)$. We count the lifts of $y_a$ with the correct invariant, and compare them with the lifts actually visited by the fine orbit.

Every lift of $y_a$ is uniquely $y_{a+1}+p^at$, where $t=(t_1,t_2)$ runs over $(\Z/p\Z)^2$. Since $Q(y)=Q(x)$, the lift has the required value of $Q$ modulo $p^{a+1}$ exactly when
\begin{equation}\label{eq:lift}
 (-2y_1-y_2)t_1+(2y_2-y_1)t_2=0\pmod p.
\end{equation}
Indeed, expand $Q(y+p^at)-Q(y)$ using \eqref{eq:diff}, divide by $p^a$, and discard the quadratic term $p^aQ(t)$ modulo $p$. The coefficient vector in \eqref{eq:lift} equals
\[
 \begin{pmatrix}-2&-1\\-1&2\end{pmatrix}y.
\]
The determinant is $-5$, prime to $p$. Also $y\equiv A^nx\not\equiv0\pmod p$, so at least one coefficient is nonzero. Choosing one coordinate freely and solving for the other gives exactly $p$ admissible lifts.

By \eqref{eq:order}, the fine period is $pK$ because $a\ge s_p$. Consider the $p$ points
\[
 A^{n+jK}x_{a+1},\qquad 0\le j<p.
\]
Each reduces to $y_a$ and preserves $Q$. They are distinct: if the points for $i$ and $j$ coincide, then the fine period $pK$ divides $(i-j)K$. Thus $p\mid i-j$, forcing $i=j$. These $p$ orbit points belong to the set of $p$ admissible lifts, so the two sets are equal. The chosen target $y_{a+1}$ is one of them, and is attained for a unique $j$. The time $n+jK$ is in the first fine circuit, as claimed.
\end{proof}
\begin{corollary}\label{cor:stable}
Fix a finite observation window $V$ of rational displacements with denominators prime to $p$, and a starting point $x$ primitive modulo $p$. The set of visited displacements stabilises as $a$ increases. Its eventual members are precisely those $v\in V$ for which $Q(x+v)=Q(x)$ over the rationals and $x_{s_p},x_{s_p}+v_{s_p}$ lie on the same orbit.
\end{corollary}
\begin{proof}
Put $\Delta_v=Q(x+v)-Q(x)$. If $\Delta_v=0$, Theorem~\ref{thm:lift} gives the same orbit-membership status at every level $a\ge s_p$. If $\Delta_v\ne0$, invariance excludes the visit for all $a>\nu_p(\Delta_v)$. Since $V$ is finite, one level suffices for all these exclusions.
\end{proof}
Thus the same integer $s_p$ marks both the end of the period plateau and the onset of stable orbit membership for pairs of equal invariant. The visited neighbours may stay the same while their directed lags change.

\subsection{Wall--Sun--Sun primes}
A prime $p>5$ is called a Wall--Sun--Sun prime if $p^2\mid F_{p-(5/p)}$, or equivalently $s_p\ge2$ \cite{Sun,Wall}. The first period plateau can be detected either by a single directed near return or by the global recovery probability.
\begin{theorem}\label{thm:wss}
Let $p>5$ be prime, $e_1=(1,0)$, $V=\{-1,0,1\}^2$, and let $P_a=\Prob(D<\infty)$ on side $p^a$. The following statements are equivalent:
\begin{enumerate}
\item $p$ is a Wall--Sun--Sun prime.
\item $K_p(1)=K_p(2)$.
\item The directed lag from $e_1$ to $(2,-1)$ is the same on sides $p$ and $p^2$.
\item The one-circuit numerator $H_V(e_1;z)$ is identical on the two grids.
\item $P_2/P_1=p^{-2}$.
\end{enumerate}
\end{theorem}
\begin{proof}
By \eqref{eq:order}, $K_p(1)=K_p(2)$ holds exactly when $s_p\ge2$, so (1) and (2) are equivalent.

For (3) and (4), the invariant at $e_1$ is $-1$. The following table gives $-Q$ on its nine observation points:
\[
\begin{array}{c|rrr}
 x_2\backslash x_1&0&1&2\\\hline
 1&-1&1&5\\
 0&0&1&4\\
 -1&-1&-1&1
\end{array}
\]
Since $p>5$, exactly three points have $Q\equiv-1\pmod{p^a}$: $e_1$, $Ae_1=(1,1)$ and $A^{-1}e_1=(2,-1)$. These three lie on the same orbit. By Lemma~\ref{lem:arithmetic}, the period is $K_p(a)$. Hence the directed lag to $A^{-1}e_1$ is $K_p(a)-1$, and by Proposition~\ref{prop:land},
\[
 H_V(e_1;z)=1+z+z^{K_p(a)-1}.
\]
The three points are distinct, so the exponent $K_p(a)-1$ is the degree of this polynomial. Either the lag or the polynomial agrees at the first two levels exactly when their periods agree. This proves the equivalence of (2), (3) and (4).

Finally, by \eqref{eq:recover},
\[
 \frac{P_{a+1}}{P_a}=\frac{K_p(a+1)}{p^2K_p(a)}
 =\begin{cases}p^{-2},&a<s_p,\\p^{-1},&a\ge s_p.\end{cases}
\]
Taking $a=1$ proves the equivalence of (1) and (5).
\end{proof}
The same plateau is therefore visible both in the lag of one particular edge and in the average probability of recovery. No Wall--Sun--Sun prime is known in the literature cited here \cite{BR}.

\section{Further directions}
For any two-by-two integer matrix $B$ of determinant one, the expression $\det(x,Bx)$ is invariant under $x\mapsto Bx$. For the cat map it equals $-Q(x)$. This suggests studying which other matrices and observation windows have near-return components of uniformly bounded size as the grid is refined.

\smallskip\noindent\textbf{Acknowledgements.} The author used Generative AI to assist with writing, programming and structural refinements.

\end{document}